\documentclass{article}
\usepackage{graphicx,color}
\usepackage{multicol,amsmath,amssymb,latexsym,amsfonts}
\usepackage{epstopdf}
\newcommand{\x}{{\bf x}}
\newcommand{\cc}{{\bf c}}
\newcommand{\n}{{\bf n}}

\newcommand{\N}{{\bf N}}

\newcommand{\s}{{\bf s}}
\newcommand{\er}{{\bf e_r}}
\newcommand{\ethe}{{\bf e_\theta}}
\newcommand{\eph}{{\bf e_\phi}}

\newcommand{\xibar}{{\overline{\xi}}}

\renewcommand{\S}{\mathcal{S}}
\newtheorem{lemma}{Lemma}
\newcommand{\pderiv}[2]{\frac{\partial #1}{\partial #2}}
\newcommand{\pderivtwo}[2]{\frac{\partial^{2} #1}{\partial #2 ^{2}}}

\newcommand{\eqr}[1]{~(\ref{#1})}
\newcommand{\figr}[1]{figure~\ref{#1}}

\newcommand{\Figr}[1]{Figure~\ref{#1}}

\newcommand{\tabr}[1]{table~\ref{#1}}

\begin{document}

\title{Fast integral equation methods for the {L}aplace-{B}eltrami equation on the sphere
\thanks{Supported in part by grants from the Natural Sciences and Engineering Research Council of Canada. NN gratefully acknowledges support from the Canada Research Chairs Council, Canada.}
}


\author{Mary Catherine A. Kropinski         \and
        Nilima Nigam 
}


%

\maketitle

\begin{abstract}
Integral equation methods for solving the Laplace-Beltrami equation on the unit sphere in the presence of multiple ``islands'' are presented. The surface of the sphere is first mapped to a multiply-connected region in the complex plane via a stereographic projection. After discretizing the integral equation, the resulting dense linear system is solved iteratively using the fast multipole method for the 2D Coulomb potential in order to calculate the matrix-vector products. This numerical scheme requires only $O(N)$ operations, where $N$ is the number of nodes in the discretization of the boundary.  The performance of the method is demonstrated on several examples. 
\end{abstract}

\section{Background}

Partial differential equations (PDEs) on surfaces and manifolds arise in many applications including image processing, biology, oceanography and fluid dynamics (see, eg., \cite{Witkin,Myers,Chaplain}). More complex topologies and metrics arise in branches of theoretical physics; \cite{lindblom}. Since solutions of these PDEs depend both on local  as well as global properties of the given differential operator on the manifold, standard numerical discretization methods (developed for PDEs in the plane or in  $\mathbb{R}^3$) need to be modified. There has been a lot of recent work in this direction, including the closest point method \cite{ruuth}, the use of surface parametrization \cite{floater}, and the use of embedding functions, \cite{Bertalmio}. Another popular strategy is to approximate the manifolds by tesselations of simpler, non-curved domains (such as triangles), and solve the PDE by projection onto these triangles,\cite{lindblom}. 

Boundary value problems on subsurfaces of the sphere arise in a variety of contexts. A by-now classical example is that of the motion of  point vortices in incompressible fluids on spheres. The fluid away from the vortex is irrotational, and is constrained  in bounded subsurfaces of the sphere with walls. Kidambi and Newton \cite{Kidambi}  used  the method of images to solve such a problem. Crowdy used stereographic projection of the subsurface into the complex plane, followed by a conformal map to the upper half plane, for the same problem \cite{Crowdy2006}. Both strategies are applicable for a rather specific set of subdomains, since they rely on the use of images and/or the knowledge of conformal maps to simpler geometries.

Just as a reduction in dimension is possible for elliptic PDEs in $\mathbb{R}^2$ or $\mathbb{R}^3$  if one reformulates the problem as an integral equation, one may reduce a boundary value problem on a subsurface of a manifold to a boundary integral equation on the boundary of the subsurface. Some prior work in this direction for the Laplace-Beltrami operator on the surface $\cal{S}$ of the unit sphere was presented in \cite{Gemmrich:2008,Gemmrich:2008p193}.

In the current paper, we present a  fast integral equation strategy for solving boundary-value problems for the Laplace-Beltrami operator associated with $\cal{S}$ in the presence of multiple ``islands''. Dirichlet boundary data is prescribed along the curves bounding these islands. Following \cite{Gemmrich:2008p193}, the boundary value problem is recast as a boundary integral equation (BIE). This reformulation is not unique: one may use direct or indirect approaches, and work with integral equations of the first or second kind. In \cite{Gemmrich:2008p193,Gemmrich:2008} numerical experiments based on Galerkin discretizations of BIE reformulations of the first kind were presented.

The analysis of of integral operators  on manifolds, or (equivalently) for elliptic problems with variable coefficients, is of course well-established, (e.g. \cite{monvel,mitrea00,mistrenistor,hsiaobook}. For example, in \cite{mitrea00} and related works,  integral equation methods are used to solve potential problems on Riemannian manifolds. However the focus is on the analytical properties of, rather than on numerical approximation strategies for, the associated boundary integral equations. Independently, in the process of solving diffraction problems from conical singularities, the authors in \cite{bonner} used an integral equation of the second kind to solve a boundary value problem on a subsurface of a sphere, but did not present an acceleration strategy.

For concreteness, we shall work with an indirect formulation based on a double layer ansatz, and examine a Fredholm integral equation of the second kind. At this juncture, we have two choices: derive an acceleration strategy for this problem on the sphere (i.e., the points on the curves live in $\mathbb{R}^3$), or use existing acceleration strategies in the complex, stereographic plane. We opt for the latter choice in this paper. 
We will show that mapping  to the stereographic plane results in an integral equation that very closely resembles the one for solving analogous boundary value problems for  Laplace's equation in the plane. 
In fact, the kernel is identical save for a ``global'' correction factor which implicitly contains information about the original manifold on which the problem is posed. This is a common feature of boundary integral operators for Laplace-Beltrami problems on manifolds which are topologically equivalent to the sphere; under stereographic projection, one obtains integral operators with the same singularities as those for the ordinary planar Laplacian, and correction terms containing information about the manifold.
Following the procedures outlined in \cite{greengard:laplace}, then, we are able to use the fast multipole method (FMM) for the two dimensional Coulomb potential \cite{THESIS,CGR,GR} to accelerate the numerical algorithm.
Such an approach follows the general solution strategy developed for a variety of linear, elliptic operators (see \cite{nishimura} for a review of FMM accelerated integral equation methods).
This strategy involves: i) formulating well-conditioned Fredholm integral equations of the second kind, ii) discretizing these integral equations with high-order quadrature schemes, iii) solving the resulting linear systems using GMRES, and iv) using FMM to compute the matrix-vector products in the iterative solution procedure.   
The result is a solution procedure that requires only $O(N)$ operations, where $N$ is the number of unknowns in the system. The method is highly accurate and is able to handle complex geometry with relative ease. 

\subsection{Generalized fundamental solution}
For the integral equation reformulation, we need  the `fundamental solution' of the Laplace-Beltrami operator $\Delta_\S$ (we recall the definition of the operator in Section 1.1). Returning to the example of vortex motion on a sphere, the 'fundamental solution' can be interpreted as the stream function of a point vortex of unit strength moving in an incompressible, irrotational fluid on the sphere without boundaries. Due to the Gauss constraint (the global integral of the vorticity must vanish), it becomes clear that we cannot have a 'fundamental solution'.  Instead, we get a   {\it generalized fundamental solution} \cite{hsiaobook} (some authors refer to this  as a {\it parametrix} \cite{firey,duduchava}), $G(\x_0,\x)$, which satisfies  the PDE $$ -\Delta_{\cal{S}} G(\x,\x_0) = \delta(\|\x-\x_0\|) -\frac{1}{4\pi}, \qquad \x\in \cal{S},
\label{eqn:lap_G}$$
and the Gauss condition for the vorticity
$$  \int_\S \Delta_\S G(\x,\x_0) ds_x =0.$$ In other words, there is a "sea" of uniform vorticity, $\frac{1}{4\pi}$, in which the point vortex at $\x_0 \in {\cal S}$ must be embedded. A generalized fundamental solution of this operator is
\begin{equation}
   G(\x,\x_0) = -\frac{1}{2\pi} \log \|\x - \x_0\| + \frac{1}{4\pi} \log 2,
   \label{eqn:greens_func}
\end{equation}
where $\|\x-\x_0\|$ is the Euclidean distance between $\x,\x_0 \in \S \subset\mathbb{R}^3$, and {\it not} the metric distance along the sphere. The parametrix $G$ satisfies\eqr{eqn:lap_G}.

\subsection{Notation and preliminary material.}
Before proceeding further, we fix some notation. The differential operators on ${\cal S}$, including the Laplace-Beltrami operator $\Delta_\S$, are most simply expressed through spherical coordinates.  A point $\x \in {\cal S}$  can be described by the spherical angles,
\[
 \x(\phi,\theta) = \left[ \begin{array}{c}
                                       \cos\phi \sin\theta \\
                                       \sin\phi \sin\theta \\
                                       \cos\theta
                                  \end{array}
                         \right] , \qquad \phi \in [0, 2\pi), \; \theta \in [0,\pi].
\]   
Let $\ethe$, and $\eph$ be the unit vectors in the $\theta$ and $\phi$ directions, respectively, and let $\er=\ethe\times\eph$. The surface gradient of a scalar field $f$ on ${\cal S}$ is
\begin{equation}
\nabla_s f(\x)  = \frac{1}{\sin\theta} \pderiv{f}{\phi} \eph + \pderiv{f}{\theta} \ethe , \label{surfgrad}
\end{equation}
The Laplace-Beltrami operator $\Delta_\S= \mbox{div}_\S \, \nabla_\S $, then, is 
 \[
 \Delta_\S  =
     \frac{1}{\sin^2 \theta} \pderivtwo{\,}{\phi} 
             + \frac{1}{\sin\theta} \pderiv{\,}{\theta} \left( \sin\theta \pderiv{\,}{\theta} \right) .
 \]

\subsection{Model problem of interest}
\begin{figure}[t]
     \centering
$\begin{array}{c}
\includegraphics[height=2.5in]{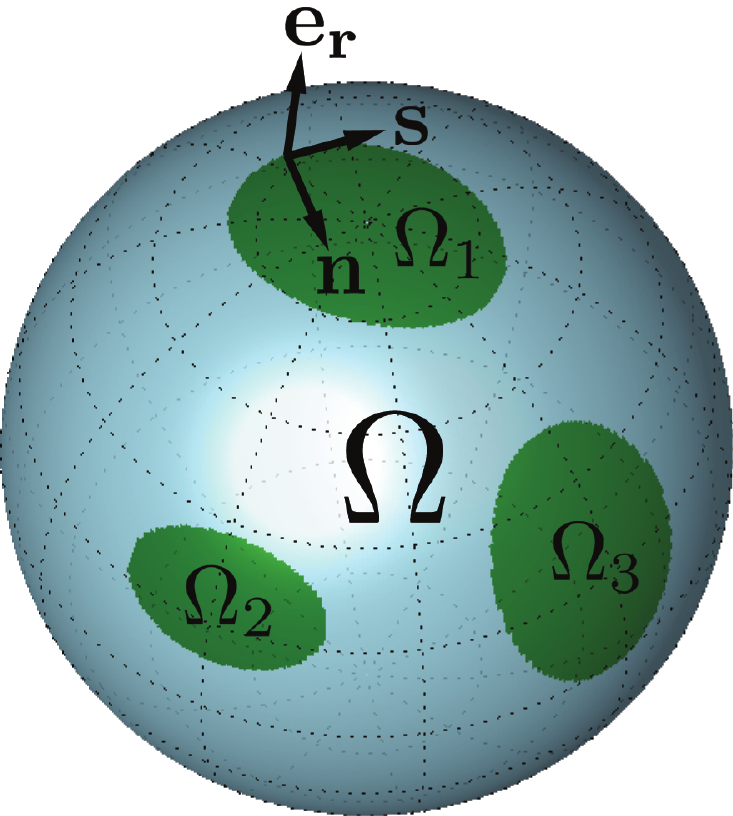}
\end{array}$
  \caption{\em The unit sphere ${\cal S}$ with $M$ ``islands'' $\Omega_1, \cdots, \Omega_M$ each bounded by smooth contours $C_1, \cdots, C_M$, respectively. On each $C_k$, $\s$ is the unit vector points in a clockwise direction with respect to $\er$, and $\n$ is the normal pointing in to $\Omega_k$ and which is also tangent to ${\cal S}$. }
  \label{fig1}
\end{figure}

In what follows, we consider the Dirichlet problem for the Laplace-Beltrami operator in multiply-connected subdomains of $\S$.  These subdomains are defined mathematically as simply-connected, sub-manifolds $\Omega_1, \cdots, \Omega_M$ of $\S$, each enclosed by a smooth boundary curve $C_1, \cdots C_M$, respectively.
Let $\Omega = {\cal S} \backslash (\Omega_1\cup \cdots \cup \Omega_M)$ be a sub-manifold on ${\cal S}$ with multiply-connected boundary curve $\partial\Omega = C_1 \cup \cdots \cup C_M$ (see \figr{fig1}). 

The Dirichlet boundary-value problem in $\Omega$ becomes: {\it find $\psi \in C^2(\bar{\Omega})$ such that} 
\begin{equation}
    \begin{array} {rll}
        \Delta_S \, \psi(\x) & = 0, & \x \in \Omega, \\
        \psi(\x) & = g(\x), & \x \in \partial\Omega.
    \end{array}
    \label{eqn:dir_bvp1}
\end{equation}

On each contour, $C_k$, oriented as shown in \figr{fig1}, the unit tangent vector is given by
\[
   {\bf s} \, ds = d\theta \, \ethe + \sin \theta \, d\phi \, \eph .
\]
The unit vector, $\n$, normal to $C_k$ and tangent to ${\cal S}$ is given by $\n = \er \times {\bf s}$, or
\begin{equation}
   \n \, ds =  d\theta \, \eph - \sin\theta d\phi \, \ethe . \label{normal}
\end{equation}

\section{Layer potentials and representation formula on the sphere}
In this section, we recall the representation formula and layer potentials for the Laplace-Beltrami equation on the sphere.
 If $\psi,\phi$ are smooth functions, and $\partial \Omega$ is a $C^1$ curve, we have the identity:
\begin{align}
   \iint_{\Omega}  & \left[ \psi(\x')  \, \Delta'_{\S} \phi(\x') - \phi(\x') \, \Delta'_{\S} \psi(\x') \right] dS' \nonumber \\
     & = \int_{\partial\Omega} \left[\psi(\x') \, \nabla'_{\S} \phi(\x') - \phi(\x') \, \nabla'_{\S} \psi(\x') \right] \cdot \n' \, ds' .   \label{green2}
\end{align}

Here, $\n'$ is the normal at $\x' \in \partial\Omega$, pointing out of the domain $\Omega$ and lying tangent to ${\cal S}$, and the primes mean differentiation or integration with respect to $\x'$.
If $\psi(\x)$ satisfies $\Delta_\S \psi=0$ on  $\Omega$,  we can use the definition of $G(\x,\x')$ from \eqref{eqn:greens_func} and \eqr{eqn:lap_G} in \eqref{green2}  to get (after standard limiting arguments) 
\[
   \psi(\x)  = \frac{1}{4\pi} \iint_{\Omega} \psi(\x') \, dS' - \int_{\partial\Omega} \left(\psi\, \nabla'_{S} G - G \, \nabla'_{S} \psi \right) \cdot \n' \, ds' . \nonumber
\]

Motivated by the integral representation in\eqr{green2}, we define  two layer potentials for the Laplace-Beltrami operator in direct analogy with layer potentials for the Laplacian in the plane, \cite{Gemmrich:2008,duduchava,mitrea00}.
\begin{description}
\item The {\it single layer potential} with sufficiently smooth density function $\rho$ is given by:
\[
   (V\rho)(\x) := \int _{\partial\Omega} \rho(\x') G(\x;\x') \, ds' .
\]
\item The {\it double layer potential} with sufficiently smooth density function $\sigma$ is given by:
\begin{align*}
  (W\sigma)(\x) & := - \int_{\partial\Omega} \sigma(\x') \, \nabla_{S'} G(\x;\x') \cdot \n' \, ds' , \\
                        & = \frac{1}{2\pi} \int_{\partial\Omega} \sigma(\x') \, 
                                     \frac{\partial \,}{\partial n'} \log ||\x - \x' || \, ds'. 
\end{align*}
\end{description}
We note that the double layer potential satisfies  the Laplace-Beltrami equation. The single layer potential $V(\rho)(\x)$ satisfies $\Delta_\S V(\rho)(\x) =\int_{\partial\Omega} \rho \, ds$; for the single layer potential to solve the Laplace-Beltrami equation, therefore, the density $\rho$ must satisfy an additional constraint that  $\int_{\partial\Omega} \rho \, ds = 0$ \cite{Gemmrich:2008p193}.
It is interesting to note, also, that the double layer potential for Laplace Beltrami on the sphere has the same form as the double layer potential for Laplace's equation in the plane. 

The basis for our boundary integral equation reformulation will be an ansatz based on the double layer potential. Two important properties of this potential are the jump relations and the behaviour of the kernel.  On smooth curves, which is the only case we report on in this paper, the double layer potential satisfies the following jump relation (see, for example, \cite{Gemmrich:2008p193} for a proof):
\begin{equation}
    \lim_{\tiny{\begin{array}{l}
        \x'  \rightarrow \x  \\
        \x'  \in  \Omega
      \end{array} }}
      (W\sigma) (\x') = \frac{1}{2} \sigma(\x) + \frac{1}{2\pi} \oint_{\partial\Omega} \sigma(\x') \, 
                                     \frac{\partial \,}{\partial n'} \log ||\x - \x' || \, ds',
      \label{eqn:jump}
\end{equation}
where $\oint$ indicates a principal-value integral.
 
We provide a simple proof that the kernel of the double layer potential is continuous in the Appendix.
\begin{lemma} On a smooth contour $C_k$, the kernel of the double layer potential is continuous. 
\end{lemma}

\begin{figure}[t]
     \centering
$\begin{array}{c}
\includegraphics[height=3.5in]{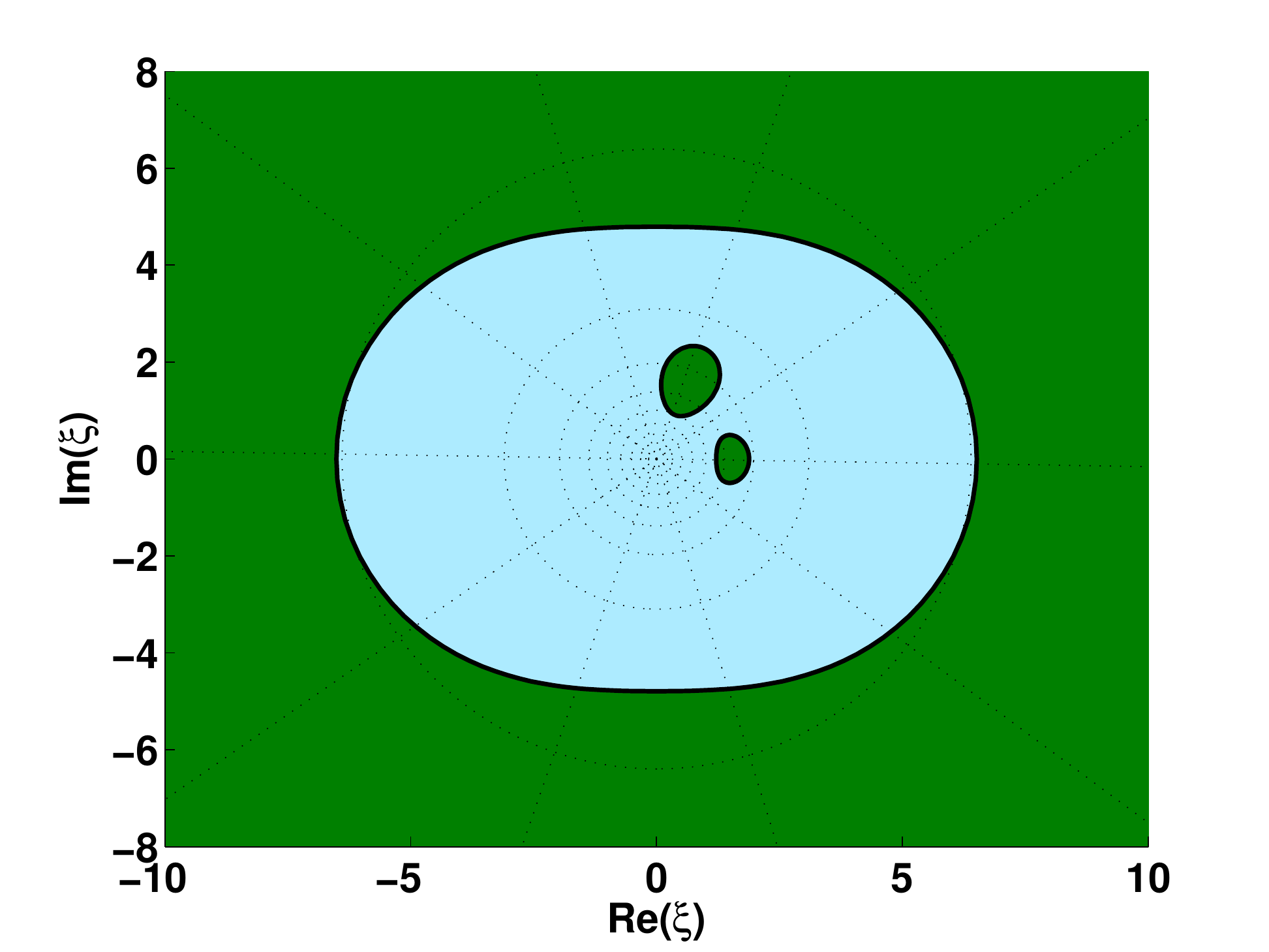}
\end{array}$
  \caption{\em The stereographic projection of the unit sphere in Figure 1. }
  \label{fig2}
\end{figure}
\section{The Stereographic Projection}
The surface of the sphere can be mapped to the complex plane through a stereographic projection. 
For $\x = (x_1,x_2,x_3)$, this mapping is given by:
\[
\xi = \cot \left( \frac{\theta}{2} \right) e^{i\phi} = \frac{x_1 + i x_2}{1 - x_3} . 
\] 
The mapping from a point $\xi$ the complex plane back to the sphere is given by
\[
   x_1 = \frac{\xi + \xibar}{1 + |\xi|^2}, \qquad 
   x_2 = \frac{\xi - \xibar}{i(1 + |\xi|^2)}, \qquad 
   x_3 = \frac{1- |\xi|^2}{1 + |\xi|^2}.
\] 

Under the stereographic projection, the North pole ($\x = (0,0,1)$) maps to the point at infinity in the complex plane, the South pole maps to the origin and the equator maps to the unit circle. 
Each of the contours, $C_k \in \S$ get mapped to contours ${\mathcal C}_k \in \mathbb{C}$. We assume without loss of generality, that the contour $C_1$ encircles the north pole. 
Hence, $\Omega$ is mapped to a region in $\mathbb{C}$ bounded by ${\mathcal C}_1$ with $M-1$ interior contours ${\mathcal C}_2, \cdots, \mathcal{C_M}$. Let $\widetilde{\Omega}$ be the mapping of $\Omega$ through the stereographic projection, let $\widetilde{\Omega_k}$ be the co\textcolor{black}{rresponding mapping of $\Omega_k$, and $\widetilde{\partial\Omega} = \mathcal{C}_1 \cup \cdots \cup \mathcal{C}_M$. }

In the stereographic plane, the Laplace-Beltrami operator is associated with the elliptic (variable coefficient) operator)
\[
   \Delta_s \equiv (1+|\xi|^2)^2 \frac{\partial^2 \,}{\textcolor{black}{\partial_{\xi\bar{\xi}}}} .
\]
\textcolor{black}{The} generalized fundamental solution for this operator in the stereographic plane \textcolor{black}{is defined} as ${\cal G}(\xi,\xi') = G(\x(\xi),\x'(\xi'))$.
\textcolor{black}{Explicitly, this distribution} becomes \cite{Crowdy2006}:
\begin{equation}
{\cal G} (\xi,\xi') = -\frac{1}{4\pi} \log \left( 2 \frac{(\xi-\xi')(\xibar-\xibar')}{(1+|\xi|^2)(1+|\xi'|^2)}
                                               \right) . \label{fundC}
\end{equation}
We proceed now with \textcolor{black}{derivation of} the mapping of the kernel of the double layer potential to the complex plane, which is most straightforward to do through the use of spherical coordinates. Using\eqr{surfgrad} and\eqr{normal} to calculate $\nabla'_S G \cdot \n'$ in spherical coordinates gives
\begin{equation}
   \frac{\partial \,}{\partial n'} G(\x,\x') = \frac{d\theta'}{\sin\theta'} \, \pderiv{G}{\phi'} - \sin\theta' d\phi' \, \pderiv{G}{\theta'} .   \label{grad1}
\end{equation}
We note, according to \cite{Crowdy:2003p30}, that 
\begin{align}
\left. \frac{\partial \,}{\partial \theta} \right|_\phi 
   & = -\frac{\xi}{\sin\theta} \left. \frac{\partial \,}{\partial \xi} \right|_{\xibar}
         -\frac{\xibar}{\sin\theta}\left. \frac{\partial \,}{\partial \xibar} \right|_{\xi}, \label{chain1} \\
\left. \frac{\partial \,}{\partial \phi} \right|_\theta 
   & = i \xi \left. \frac{\partial \,}{\partial \xi} \right|_{\xibar}
         -i \xibar \left. \frac{\partial \,}{\partial \xibar} \right|_{\xi} . \label{chain2}
\end{align}     
These two relations similarly hold with respect to primed variables $\theta'$, $\phi'$, $\xi'$, and $\xibar'$.    
Substituting\eqr{chain1} and\eqr{chain2} into\eqr{grad1} yields
\begin{align}
 \frac{\partial \,}{\partial n'} G(\x,\x')
    = -& \frac{i\xi'}{\sin\theta'} \frac{\partial{\cal G}}{\partial \xi'} (-d\theta' + i \sin\theta' d\phi')  + \frac{i\xi'}{\sin\theta'} \frac{\partial{\cal G}}{\partial \xibar'} (-d\theta' - i \sin\theta' d\phi') .
       \label{grad2}
\end{align}
Differentiating $\xi = \cot \left( \frac{\theta}{2} \right) e^{i\phi}$ and simplifying gives
\begin{equation}
   \frac{\sin\theta d\xi}{\xi} = -d\theta + i \sin\theta d\phi . \label{dxi}
\end{equation}
Substituting\eqr{dxi} into\eqr{grad2} gives 
\begin{equation}
 \frac{\partial \,}{\partial n'} G(\x,\x')
    = - i \frac{\partial \, }{\partial \xi'} {\cal G}(\xi,\xi') d\xi' 
       + i \frac{\partial \,}{\partial \xibar'} {\cal G}(\xi,\xi') d\xibar' ,
       \label{grad3}
\end{equation}
and finally, after differentiating\eqr{fundC}, we can therefore rewrite the double layer potential in the stereographic plane as 
\begin{equation}
   (W\sigma)(\xi) = \mbox{Re} \left\{ \frac{1}{2\pi i} \int_{\partial\Omega} \sigma(\xi')  
                              \left[ \frac{1}{\xi-\xi'} - \frac{\bar{\xi'}}{1 + |\xi'|^2} \right] d\xi' \right\} .
                              \label{dlpC}
\end{equation}
The first term in the above is simply the double layer potential for Laplace's equation in the plane. 
We note here for future reference that the double layer potential is continuous for smooth curves, and that
\begin{equation}
	\lim_{\xi'\rightarrow\xi} \mbox{Im} \left\{  \frac{d\xi'}{\xi-\xi'}\right\}
	   = \frac{1}{2} \kappa(\xi) |d\xi|, \label{cont_kern}
\end{equation}
where $\kappa(\xi)$ is the curvature of $\partial\Omega$ at the point $\xi$ \cite{greengard:laplace}.

\section{The Integral Equation reformulation}
In the case of a simply-connected domain, we employ the double-layer ansatz, $\psi(\x) = (W\sigma)(\x)$. Using the jump relation\eqr{eqn:jump}, together with the boundary conditions,  we obtain the integral equation for the solution to\eqr{eqn:dir_bvp1}  as
\begin{equation}
\frac{1}{2} \sigma (\x) + \frac{1}{2\pi} \oint_{\partial\Omega} \sigma(\x') \frac{\partial\,}{\partial n'} \log\|\x-\x'\|  \, ds' = g(\x) . \label{inteqn:simple}
\end{equation}
where $\oint$ indicates a principal-value integral, or equivalently in the stereographic plane 
\begin{equation}
\frac{1}{2} \sigma (\xi ) + \mbox{Re} \left\{ \frac{1}{2\pi i} \int_{\widetilde{\partial\Omega}}
                                 \sigma(\xi') \, \left[ \frac{1}{\xi-\xi'} - \frac{\xibar'}{1+|\xi'|^2} \right] \, d\xi' \right\} 
    = g(\x(\xi)). \label{stereo:simple}
\end{equation}
The Cauchy kernel in the latter form of the integral equation makes it easier to determine the nontrivial homogeneous solutions in the multiply-connected case.
Since the kernels in the above integral equations are continuous if the boundary $\partial \widetilde{\Omega}$ is smooth, \textcolor{black}{both} \eqr{inteqn:simple} and\eqr{stereo:simple} are Fredholm equations of the second kind with \textcolor{black}{compact integral operators}.  Moreover, if $\Omega$ is simply connected these equations have no non-trivial homogeneous solutions. Therefore, from the Fredholm alternative, eqr{inteqn:simple} and\eqr{stereo:simple} have a unique solution for given integrable data $g(\x)$. 

\subsection{The Multiply-Connected Case}
There is a fundamental difficulty with using \eqr{inteqn:simple} and\eqr{stereo:simple} when $\partial\Omega$ is $(M-1)$-ply connected, namely that there are $M-1$ nontrivial homogeneous solutions. 
The nature of these nontrivial homogeneous solutions is most straightforward to see in\eqr{stereo:simple} because of the resemblance of this equation to the integral equation based on the double layer potential for Laplace's equation in the plane. 
To see this, define $M$ layer densities as  
\[
   \zeta_i (\xi) = \left\{ \begin{array}{rcl} 
                                        1 & \mbox{for} & \xi \in \mathcal{C}_i , \\
                                        0 & \mbox{for} & \xi \in \mathcal{C}_k, k\ne i, 
                                 \end{array}
                       \right.
                                 \qquad i = 1, \cdots, M. 
\]
Then, it follows directly from well-known identities that 
\[
   (W\zeta_1)(\xi) = \left\{ \begin{array}{rcl}
                       1 + D_1, & \mbox{for} & \xi \in \widetilde{\Omega}_1,\\
                       \frac{1}{2} + D_1, & \mbox{for} & \xi \in \mathcal{C}_1 ,\\
                       D_1, & &\mbox{elsewhere} , 
                 \end{array}
         \right.
\]
and for $k> 1$, 
\[
   (W\zeta_k)(\xi) = \left\{ \begin{array}{rcl}
                       -1 + D_k, & \mbox{for} & \xi \in \widetilde{\Omega}_1,\\
                       - \frac{1}{2} + D_k, & \mbox{for} & \xi \in \mathcal{C}_1 ,\\
                       D_k, & &\mbox{elsewhere} . 
                 \end{array}
         \right.
\]
In the above,
\[
D_k = \mbox{Re} \left\{ -\frac{1}{2\pi i} \int_{\widetilde{\Omega}_k}
                                 \frac{\xibar'}{1+|\xi'|^2} \, d\xi' \right\} = \frac{1}{\pi} \int_{Int(\mathcal{C}_k)}   \frac{1}{1+x^2+y^2} \, dx \, dy  .
\]
There are $M-1$ non-trivial homogeneous solutions of\eqr{inteqn:simple} in the form:
\[
   \tilde{\sigma}_i (\xi) = \left\{ \begin{array}{lcl} 
                                        1, & \mbox{for} & \xi \in \mathcal{C}_1 , \\
                                        -(1+D_1)/D_i, & \mbox{for} & \xi \in \mathcal{C}_i, \\
                                        0, & \mbox{for} & \xi \in \mathcal{C}_k, k \ne i, 
                                 \end{array}
                       \right.
                                 \qquad i = 2, \cdots, M. 
\]
Clearly $D_i\not=0$ unless the curve $\mathcal{C}_i$ is degenerate, ie, has zero enclosed area.
Thus, the solution to\eqr{inteqn:simple} is non-unique. We resolve this issue for \eqr{eqn:dir_bvp1} as follows. 

From Fredholm theory, \eqr{inteqn:simple} can be solved only if the right-hand side is orthogonal to the $M-1$ independent solutions of the adjoint integral equation. 
This would require the determination of the $M-1$ independent solutions of the adjoint system. 
Following the approach taken in \cite{greengard:laplace,Mikh}, \textcolor{black}{we instead} seek a solution in the form
\begin{equation}
\psi(\x) = \frac{1}{2\pi} \int_{\partial\Omega} \sigma(\x') \frac{\partial\,}{\partial n'} \log\|\x - \x'\|  \, ds' + \sum_{k=1}^M {\cal A}_k \, G(\x,\cc_k), \label{eqn:u_multi}
\end{equation}
where $\cc_k \in \Omega_k$. 
We can view this physically as placing a fundamental solution of unknown strength inside each $\Omega_k$. Mathematically, these singularities play the role of Lagrange multipliers.
Note that for $\x \in \Omega$, and using\eqr{eqn:lap_G}, 
\[
   \Delta_S \psi = \frac{1}{4\pi} \sum_{k=1}^M {\cal A}_k.
\]
\textcolor{black}{In order for $\psi$} to satisfy $\Delta_\S \psi = 0$, we \textcolor{black}{impose} the following constraint on the \textcolor{black}{strengths of the singularites}:
\begin{equation}
       \sum_{k=1}^M {\cal A}_k = 0 .
       \label{eqn:sum_ak}
\end{equation}
Applying the jump relations as before, the resulting integral equation is 
\begin{equation}
\frac{1}{2} \sigma (\x) + \frac{1}{2\pi} \oint_{\partial\Omega} \sigma(\x') \frac{\partial\,}{\partial n'} \log\|\x-\x'\|  \, ds' 
     + \sum_{k=1}^M {\cal A}_k \, G(\x,\cc_k)  = g(\x) . \label{inteqn:multi_1}
\end{equation}
Given the $M-1$ extra degrees of freedom (\eqr{eqn:sum_ak} must be satisfied), we augment \textcolor{black}{the integral equation \eqr{inteqn:multi_1}} with $M - 1$ constraints of the form 
\begin{equation}
\int_{C_k} \sigma(\x) \, ds = 0, \qquad k=2, \cdots, M. \label{constraints}
\end{equation}

For reasons that will be discussed in the following section, we will be discretizing the integral equation in the stereographic plane. Summarizing, the system of equations that correspond to solving\eqr{eqn:dir_bvp1} is:
\begin{equation}
\begin{array}{rcl}
\lefteqn{\frac{1}{2} \sigma (\xi ) + \mbox{Re} \left\{ \frac{1}{2\pi i} \int_{\widetilde{\partial\Omega}}
                                 \sigma(\xi') \, \left[ \frac{1}{\xi-\xi'} - \frac{\xibar'}{1+|\xi'|^2} \right] \, d\xi' \right\}    }\hspace{2.5in} \medskip \\
      + {\displaystyle \sum_{k=1}^M {\cal A}_k \, {\cal G} (\xi,\xi_k )} & = & g(\x(\xi)),   \medskip \\
{\displaystyle \int_{{\cal C}_k} \sigma \, d\alpha  }& = & 0, \qquad k = 2, \cdots, M, \medskip \\
{\displaystyle \sum_{k=1}^M {\cal A}_k } & = & 0, 
\end{array}
\label{continuous_sys}
\end{equation}
where $\xi_k$ is the mapping of ${\bf c_k}$ to the stereographic plane. 
Here, the constraint\eqr{constraints} has been replaced by integrating with respect to the curve parameter $\alpha$ (the parametrization of ${\cal C}_k$ will be discussed in detail in the next section).
By the preceding discussion, this system is invertible and simultaneously determines the values of the constants ${\cal A}_k$ and the desired layer density $\sigma(\xi)$. 
We will show in the numerical results that after suitable preconditioning, the discrete linear system corresponding to the discretization of the above has a bounded condition number.
%

\section{Numerical Methods}
\textcolor{black}{We now describe the discretization and solution of the boundary integral reformulation \eqr{continuous_sys} of \eqr{eqn:dir_bvp1}.} 
We will denote by $\sigma^k$ and $\xi^k$  the restriction of \textcolor{black}{$\sigma$ and $\xi$} to the component curve ${\cal C}_k$, $k=1,...,M$. 
We assume each component curve is parametrized with ${ \xi}^k(\alpha):[0,2\pi) \rightarrow {\cal C}_k$.
We discretize\eqr{inteqn:multi_1} using the Nystr\"{o}m 
discretization
based on the trapezoidal rule.
For this, we assume we are given $N$ points on each contour ${\cal C}_k$, equispaced with respect to $\alpha$.  This strategy will achieve super-algebraic convergence \textcolor{black}{in N}
for smooth data on smooth boundaries ${\cal C}_k$.
Associated with each such point \textcolor{black}{$\xi^k_j \in \mathcal{C}^k$, $j=1,....N, k=1,...,M$}, is an unknown density $\sigma^k_j$, and $d\xi^k_j$, the derivative of $\xi^k$ with respect to $\alpha$ at that point. 
The mesh spacing is $h=2\pi/N$ and the total number of points is $NM$.
Discretizing the integral equation in\eqr{continuous_sys} yields the following $NM$ equations:
\begin{equation}
\sigma^k_i + 2 h \sum_{m = 1}^M \sum_{j = 1}^{N} \sigma_j K(\xi^k_i,\xi^m_j) 
  +  2\sum_{m=1}^M {\cal A}_m G(\xi_i^k ,\xi_m) =  2 g^k_i, \quad k=1,...,M, i=1,...,N . \label{eqn:disc_sphere}
\end{equation}
where
\begin{align*}
   K(\xi^k_i,\xi^m_j) 
	   & = \frac{1}{\pi} \, \mbox{Im}
		\left\{ \frac{d\xi^m_j}{\xi^k_i - \xi^m_j} 
		       - \frac{\overline{\xi_j^m} \, d\xi_j^m}{1+|\xi^m_j|^2} \right\} , \; \; \xi^k_i \ne \xi^m_j\\
   K(\xi^k_i,\xi^k_i) 
	   & = \frac{1}{4\pi} |d\xi^k_j| \kappa^k_j 
		    - \frac{1}{\pi}\mbox{Im} 
			   \left( \frac{\overline{\xi^k_j} \, d\xi^k_j} {1 + |\xi^k_i|^2}      
				\right).
\end{align*}
The first term in $K(\xi^k_i,\xi^k_i)$ is derived from\eqr{cont_kern}, and $\kappa_i^k$ is the curvature at the point $\xi^k_i$
Discretizing the constraints in\eqr{continuous_sys} yields
\begin{equation}
 \sum_{j=1}^N \sigma^k_j = 0, \; \; k = 2, \cdots, M.\label{eqn:disc_constr}
\end{equation}
Equation\eqr{eqn:disc_sphere}, together with\eqr{eqn:disc_constr} and\eqr{eqn:sum_ak} give $(N+1)M$ equations for $(N+1)M$ unknowns. We can rewrite this discrete systems as
\begin{equation}
\left[ \begin{array}{cc}
          I + K & E \\
          F & D \end{array}
\right]
\left[ \begin{array}{c} {\vec \sigma} \\ {\vec a} \end{array} \right] =
\left[ \begin{array}{c} 2 {\vec g} \\ 0 \end{array} \right], \label{matrix_eqns}
\end{equation}
where 
\[
{\vec \sigma} = ({\sigma}^1_1, \cdots, {\sigma}^1_N, \cdots, 
{\sigma}^M_1, \cdots, {\sigma}^M_N)^T,
\]
is the vector of unknown layer densities, 
\[
{\vec a} = ({\cal A}_1, {\cal A}_2, \cdots, {\cal A}_M)^T,
\] 
is the vector of unknown singular source strengths, and 
\[
{\vec g} = (g_1^1, \cdots, g_N^1, \cdots, g_1^M, \cdots, g_N^M)^T
\]
is the vector of given boundary values. 
The $NM\times M$ matrix $E$ represents the influence of the constants ${\cal A}_k$ on the potential field, and the $M\times NM$ matrix $F$ represents the discrete constraint equations in\eqr{eqn:disc_constr}. The $M\times M$ block matrix $D$ incorporates the constraint\eqr{eqn:sum_ak}. 

The matrix equations\eqr{matrix_eqns} are solved iteratively using GMRES. Since the integral operator in\eqr{inteqn:multi_1} is compact, $I+K$ is a low rank perturbation of the identity. \textcolor{black}{In addition, as is typical with discretizations of integral operators, the matrices $K$ are dense $NM\times NM$ matrices.}  In order to improve the convergence rate of GMRES, we use a preconditioner that eliminates the influence of the singular source terms. We note that a very similar procedure was used in \cite{greengard:laplace} and \cite{GKM}. For completeness, we include a description of this preconditioning procedure and refer the reader to \cite{greengard:laplace} for further details.
Instead of solving\eqr{matrix_eqns}, we \textcolor{black}{iteratively} solve
\begin{equation}
\left[ \begin{array}{cc}
          I & E \\
          F & D \end{array}
\right]^{-1}
\left[ \begin{array}{cc}
          I + K & E \\
          F & D \end{array}
\right]
\left[ \begin{array}{c} {\vec \sigma} \\ {\vec a} \end{array} \right] =
\textcolor{black}{\left[ \begin{array}{cc}
          I  & E \\
          F & D \end{array}
\right]^{-1} } 
\left[ \begin{array}{c} 2 {\vec g} \\ 0 \end{array} \right]. \label{matrix_eqns_c}
\end{equation}
At each iteration, the preconditioning matrix \[
 \textcolor{black}{\left[ \begin{array}{cc}
          I  & E \\
          F & D \end{array}
\right]} 
\]
\textcolor{black}{ must be inverted. } This is accomplished by first forming the $M\times M$ Schur complement $S$ of $D$ in the preconditioner,
given by
\[
S = D-FE .
\]
The Schur complement $S$ is factored using Gaussian elimination, which requires $\frac{1}{3} M^3$ operations. To solve the linear system
\[
\left[ \begin{array}{cc}
          I  & E \\
          F & D \end{array}
\right]
\left[ \begin{array}{c} {\vec z}_\sigma \\ {\vec z}_a \end{array} \right] =
\left[ \begin{array}{c} {\vec r}_\sigma \\ {\vec r}_a \end{array} \right] 
\]
we first backsolve to find ${\vec z}_a$ from
\[
S {\vec z}_a = {\vec r}_a - F {\vec r}_\sigma,
\]
which requires $O(M^2)$ operations, and then compute the values ${\vec z}_\sigma$ from
\[
{\vec z}_\sigma = {\vec r}_\sigma - E {\vec r}_a
\]
which requires $O(NM)$ operations.
As shown in \cite{greengard:laplace}, the result of this preconditioner is that the eigenvalues of the system matrix become more tightly clustered around 1. We will demonstrate in \textcolor{black}{Example 1} that \textcolor{black}{the use of  this preconditioning strategy results in a bounded  (in $N$) condition number for the preconditioned system. Preconditioning also leads to a significant reduction in both the number of GMRES iterations and overall CPU time.}

As with solving most linear systems which arise from the discretization of an integral equation, the bulk of the computational cost comes from calculating the matrix-vector product in\eqr{matrix_eqns}.
By moving to the stereographic plane, we can exploit the two dimensional FMM based on the Coulomb potential which enables us to compute $[I+K] {\vec \sigma}$ in $O(NM)$ time. 
We refer the reader to the original papers \cite{CGR,GR,THESIS} on FMM for further details.
We could also use FMM to evaluate $E {\vec a}$; however, our intended applications will have moderate values of $M$, and so we calculate this product direcly, in $O(NM^2)$ time. 
The number of iterations required to solve\eqr{matrix_eqns} to a fixed precision $\epsilon$ is bounded, at $c(\epsilon)$, say, the total number of operations for solving the preconditioned system is approximately
\[
 \frac{1}{3} M^3 + c(\epsilon) (M^2 + M^2 N)  .
\]


\section{Numerical Results}
\textcolor{black}{We now demonstrate the effectiveness of the strategy described in the previous section. In what follows, we consider some boundary value problems on $\S$ where the exact solution can be analytically computed. This allows us to evaluate the errors of discretization explicitly. In the first example we compare the performance of an iterative solver on  \eqr{matrix_eqns}, with and without the preconditioner. In Example 2, we consider a more geometrically complex region - $\Omega$ represents the oceans away from the major continents of the Earth. In Example 3, we demonstrate the effectiveness of the acceleration technique on a problem with a large number of unknowns.}

The algorithms described above have been implemented in Fortran. The tolerance for convergence of GMRES is set to $10^{-11}$, with GMRES restarting every 50 iterations.
The FMM tolerance is set to $10^{-14}$. 
Here, we illustrate the performance on a variety of examples.
All timings cited are for a single processor on a Mac Pro 4.1 with two 2.93 GHz Quad-Core Intel Xeon processors. 

\begin{figure}
\centering
$\begin{array}{cc}
\includegraphics[width=2.in]{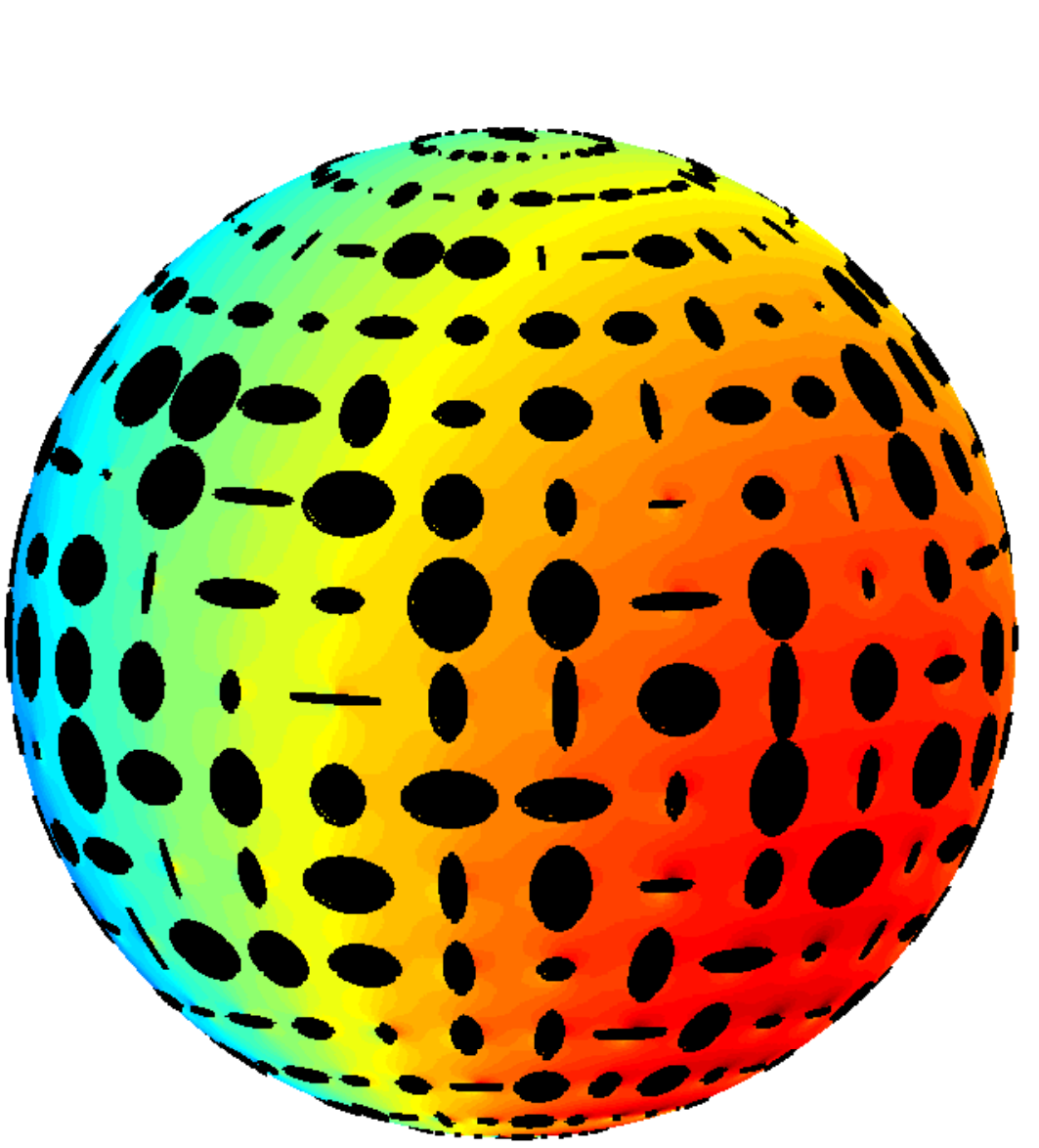} \hspace{.25in} & \hspace{.25in}
\includegraphics[width=2.in]{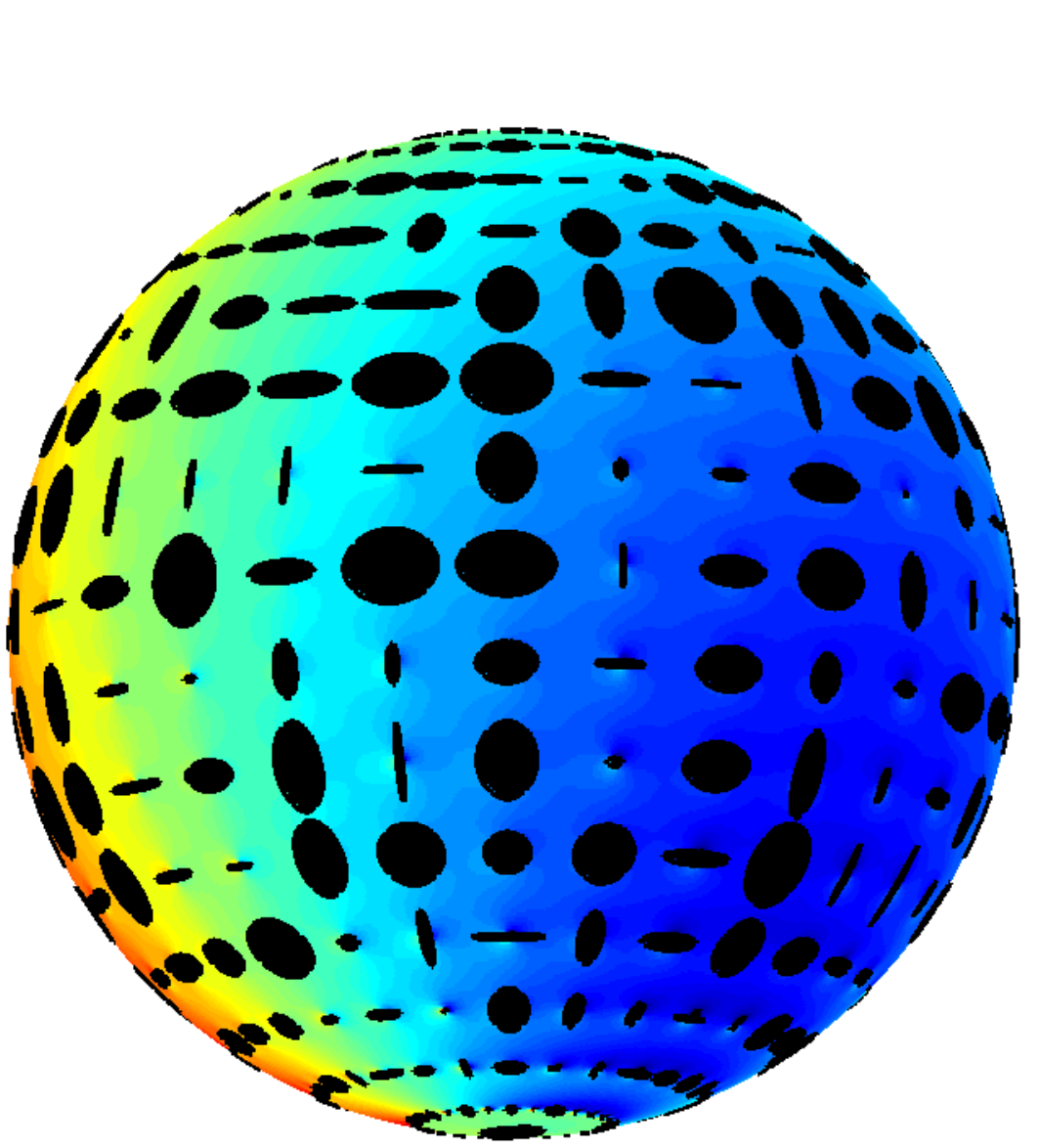} \\
\includegraphics[width=2.in]{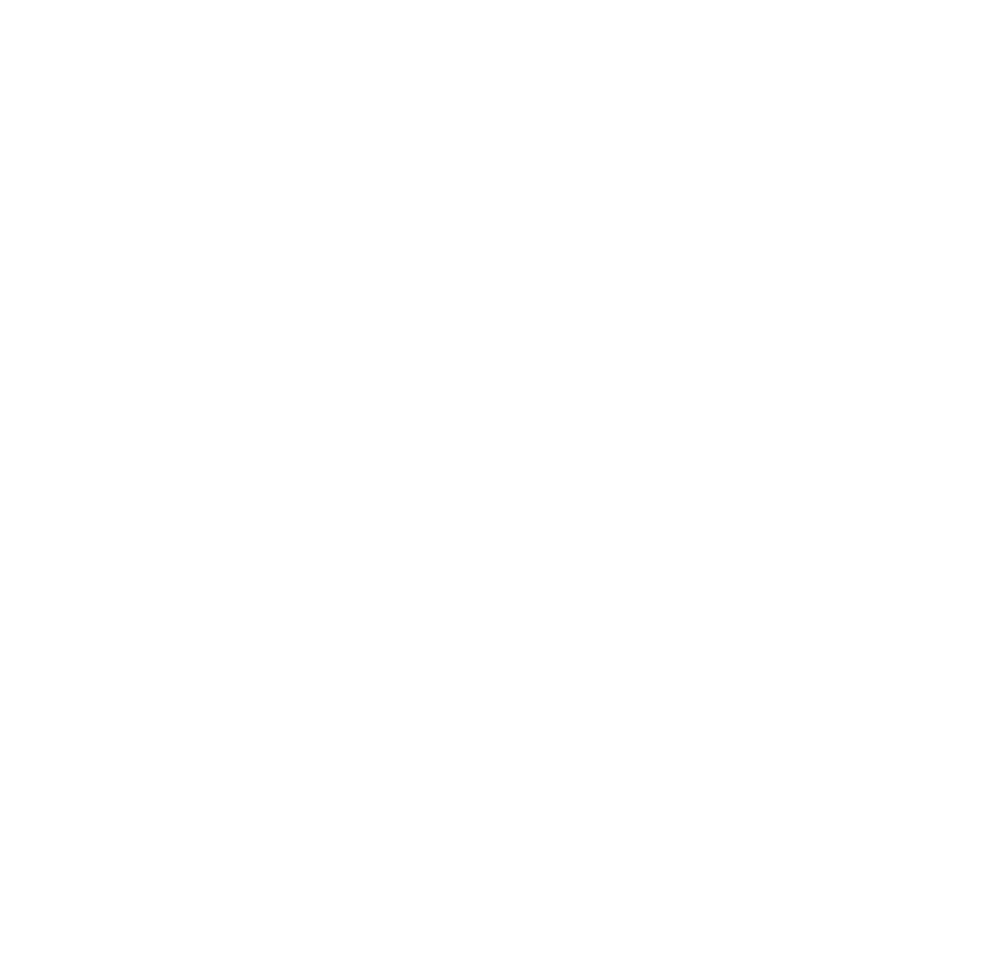} \hspace{.25in}  & \hspace{.25in}
\includegraphics[width=2.in]{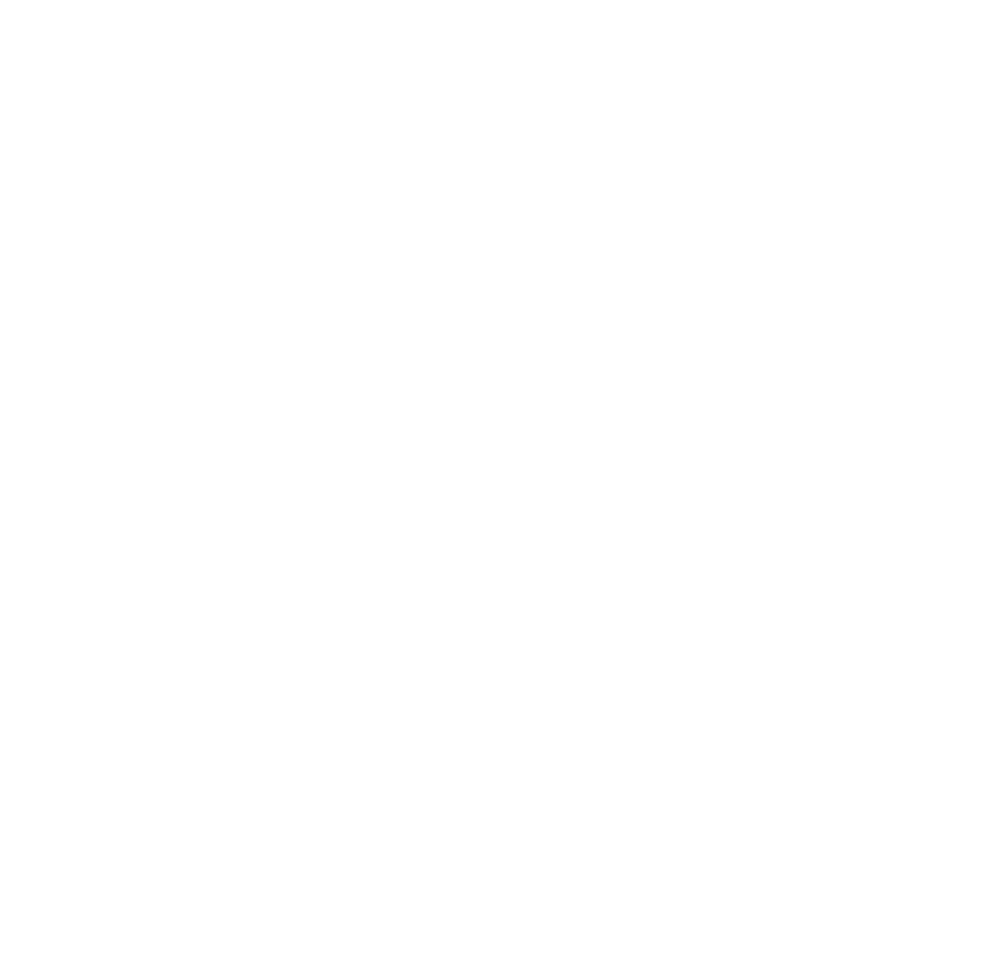}
\end{array}$
\caption{\em The solution to Example 3 with $M=407$. The top two plots show two different view points of the solution on the sphere, and the bottom shows the solution in the stereographic plane. }
\label{fig3}
\end{figure}

\begin{figure}[p]
\centering
$\begin{array}{cc}
\includegraphics[width=2.in]{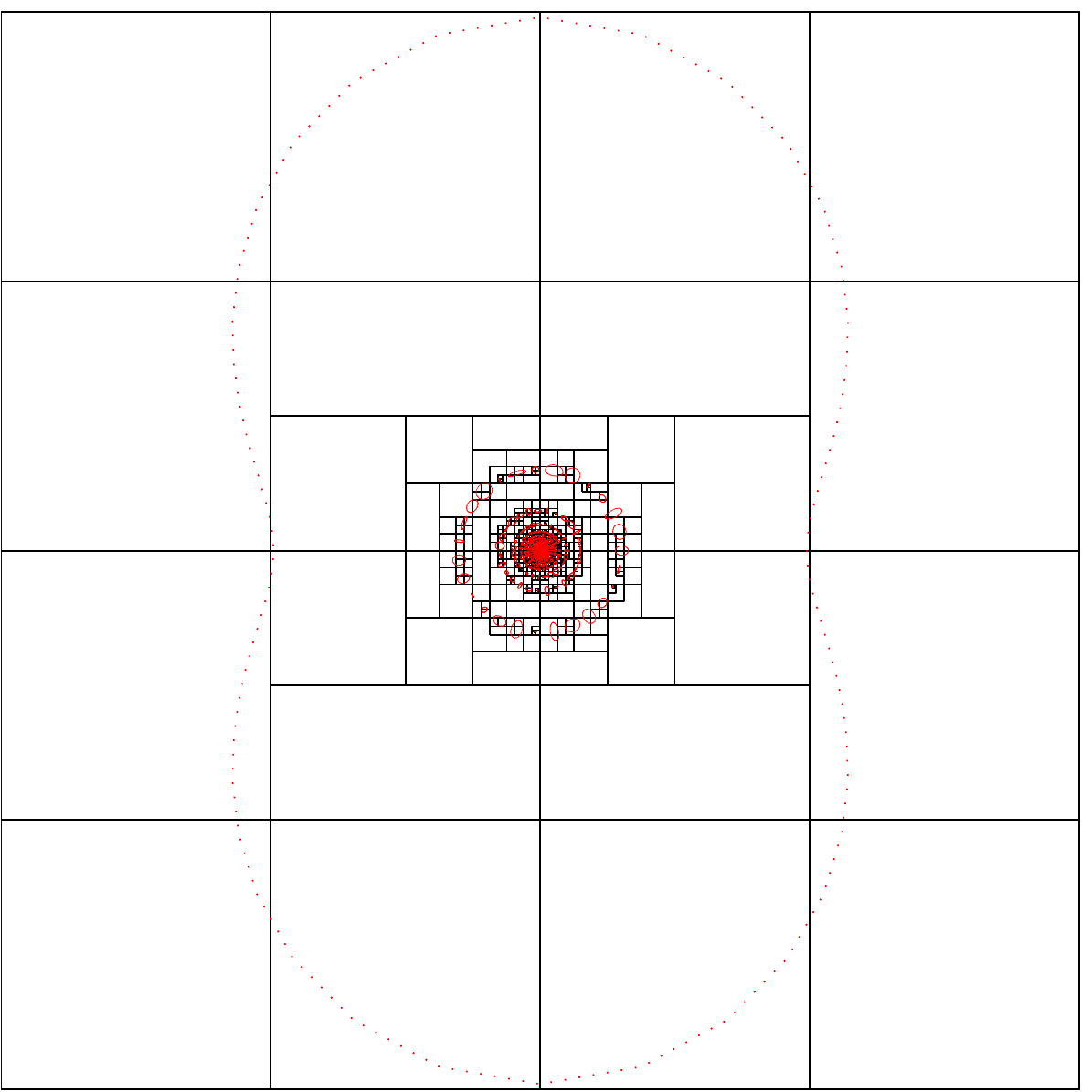} \hspace{.25in} & \hspace{.25in}
\includegraphics[width=2.in]{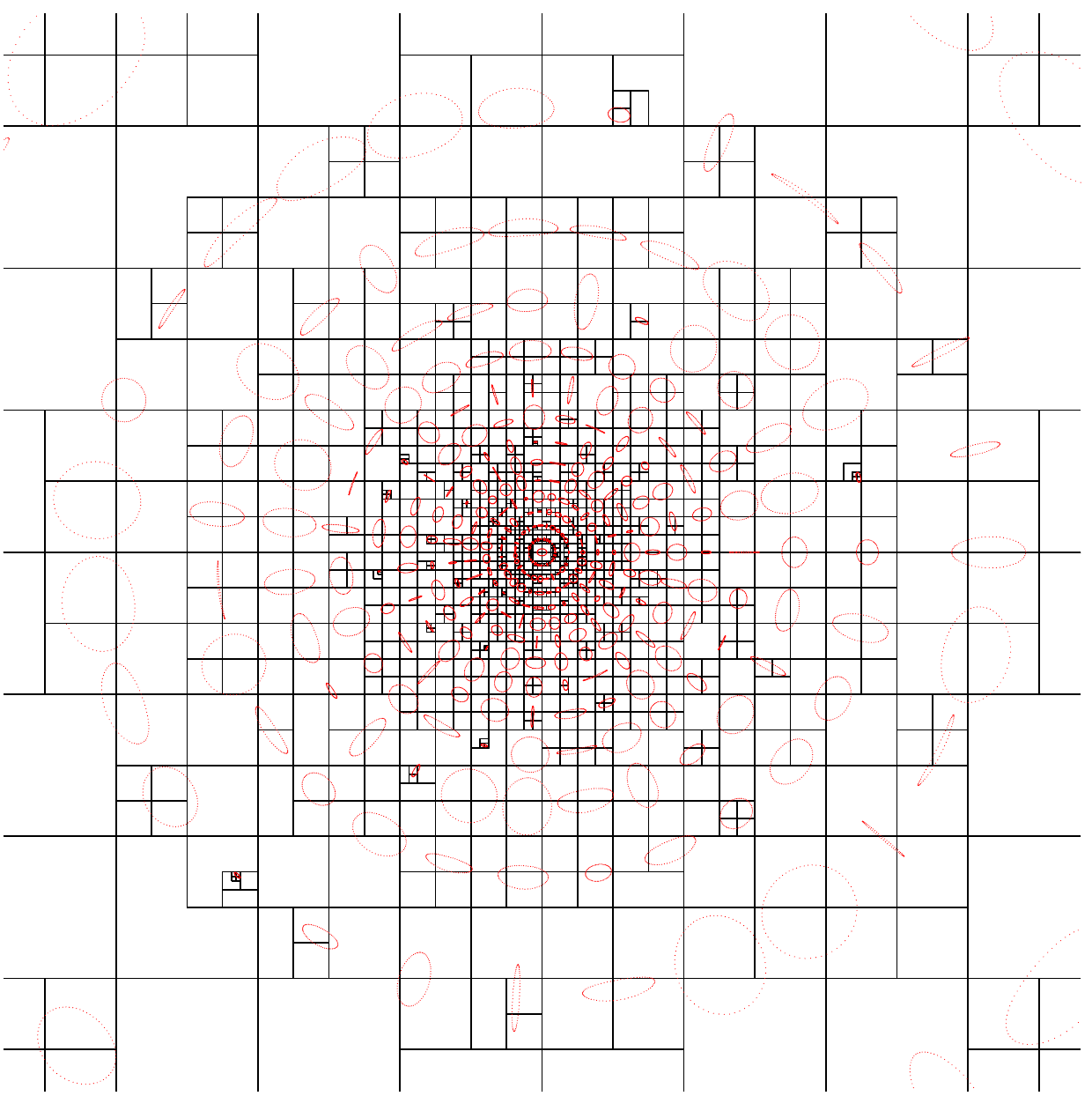} 
\end{array}$
\caption{\em The domain in the stereographic plane with $M=407$, including the quadtree structure for FMM. The plot on the left is the entire domain and the plot on the right is a close-up view near the origin (corresponding to the south pole on the sphere).}
\label{fig4}
\end{figure}

\noindent {\sc Example 1: Preconditioning.}
We first consider the problem of solving the Laplace Beltrami equation on the sphere in the presence of M elliptical islands, where the boundary conditions in the stereographic plane are generated by 
\begin{equation}
 \psi(\xi,\xibar)  = \frac{1}{2} \sum_{k=1}^M \mbox{Re} \frac{1}{\xi-\xi_k} .
 \label{exact_bcs}
 \end{equation}

To compare the performance of solving the unpreconditioned versus preconditioned system, we  consider the case with $M=15$.
The computational performance is outlined in \tabr{table1}.
As can be seen in this table, the condition number of the system matrix increases linearly \textcolor{black}{in $N$} if preconditioning is not used. With preconditioning, the conditioning \textcolor{black}{of the system} is improved significantly and the condition number remains bounded \textcolor{black}{in $N$}.

\begin{table}[b]
\caption{Performance of the algorithm on example 1 with $M=15$. The error is the maximum relative error in the solution sampled at 75 points distributed throughout the sphere, CPU is the total time taken by GMRES, $K$ is the condition number of the system matrix and  \# is the number of GMRES iterations required for convergence.
\label{table1} }
\begin{tabular*}{4in}{@{\extracolsep{\fill}}rrrrl }    
\hline\noalign{\smallskip}
      \multicolumn{5}{c}{{\bf Unpreconditioned} }  \\
\multicolumn{1}{c}{$N$} & 
\multicolumn{1}{c}{\#} &  \multicolumn{1}{c}{$K$} &\multicolumn{1}{c} {CPU} 
& \multicolumn{1}{c}{Error }  \\
\noalign{\smallskip}\hline\noalign{\smallskip}
32   & 88 & 2151.2 & 2.7 & $6.522 \times 10^{-3}$ \\ 
64   & 88 & 4230.9 & 3.0 & $3.025\times 10^{-5}$ \\  
128 & 88 & 8389.7 & 3.4 & $6.307\times 10^{-10}$ \\  
256 & 88 & 16706.8 & 5.2 & $1.828\times 10^{-11}$ \\  
512 & 88 & 33341.0 & 5.6 & $1.812\times 10^{-11}$ \\  
\noalign{\smallskip}\hline
\end{tabular*}

\vspace{0.2in}
\begin{tabular*}{4in}{@{\extracolsep{\fill}}rrrrl }  
\hline\noalign{\smallskip}
       \multicolumn{5}{c}{{\bf Preconditioned} } \\
\multicolumn{1}{c}{$N$} & 
\multicolumn{1}{c}{\#} &  \multicolumn{1}{c}{$K$} &\multicolumn{1}{c} {CPU} 
& \multicolumn{1}{c}{Error }  \\
\noalign{\smallskip}\hline\noalign{\smallskip}
32    & 26 & 35.4 & 0.9 & $6.522 \times 10^{-3}$    \\ 
64     & 26 & 35.4 & 0.9 & $3.025\times 10^{-5}$   \\  
128  & 26 & 35.3 & 1.1 & $6.612\times 10^{-10}$  \\  
256  & 26 & 35.3 & 1.4 & $1.020\times 10^{-10}$\\  
512 & 26 & 35.3 & 1.9& $1.132\times 10^{-10}$\\  
\noalign{\smallskip}\hline\end{tabular*}
\end{table}

\begin{figure}
\centering
$\begin{array}{cc}
\includegraphics[width=2.in]{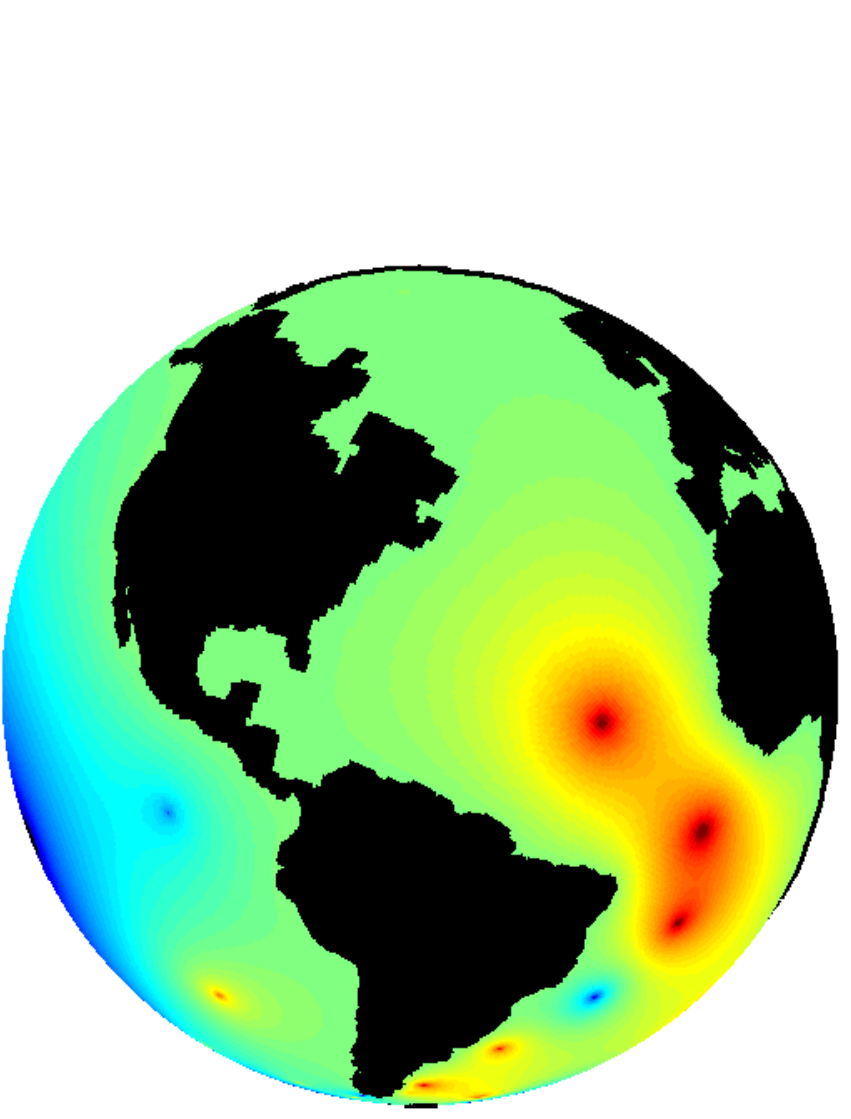} \hspace{.25in} & \hspace{.25in}
\includegraphics[width=2.in]{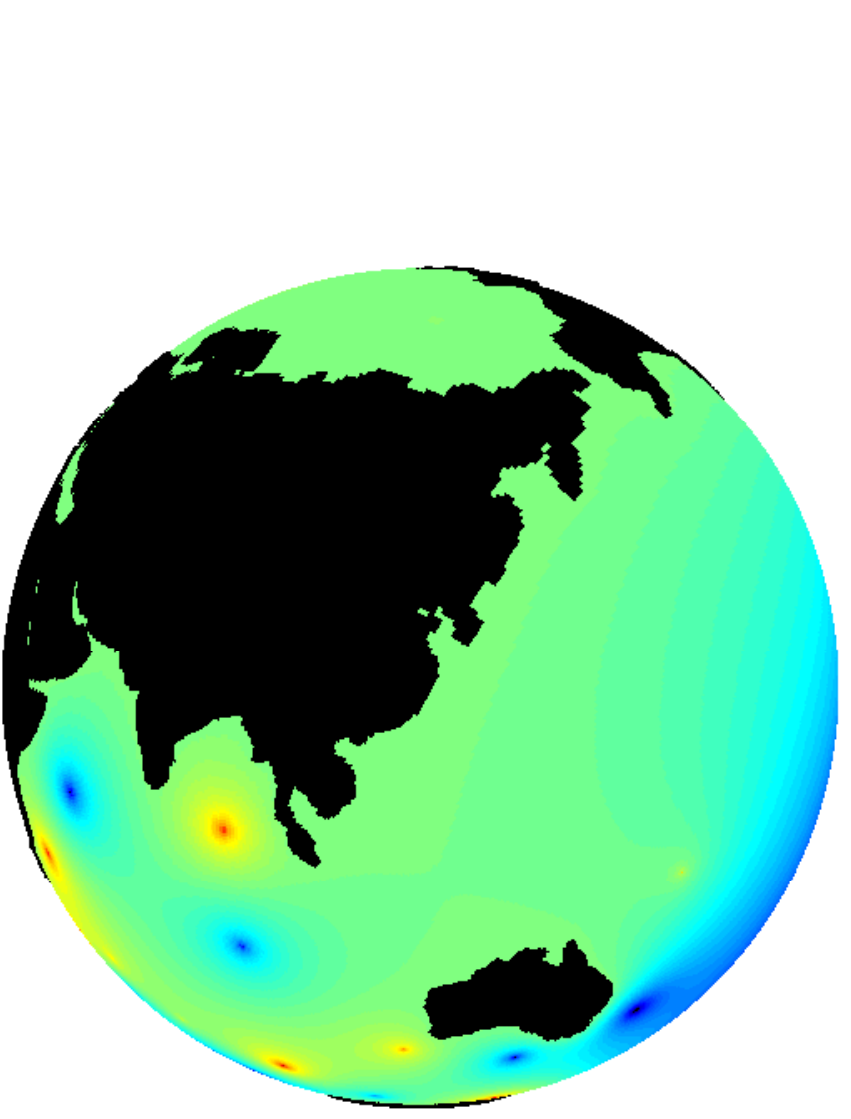} \\
\includegraphics[width=2.in]{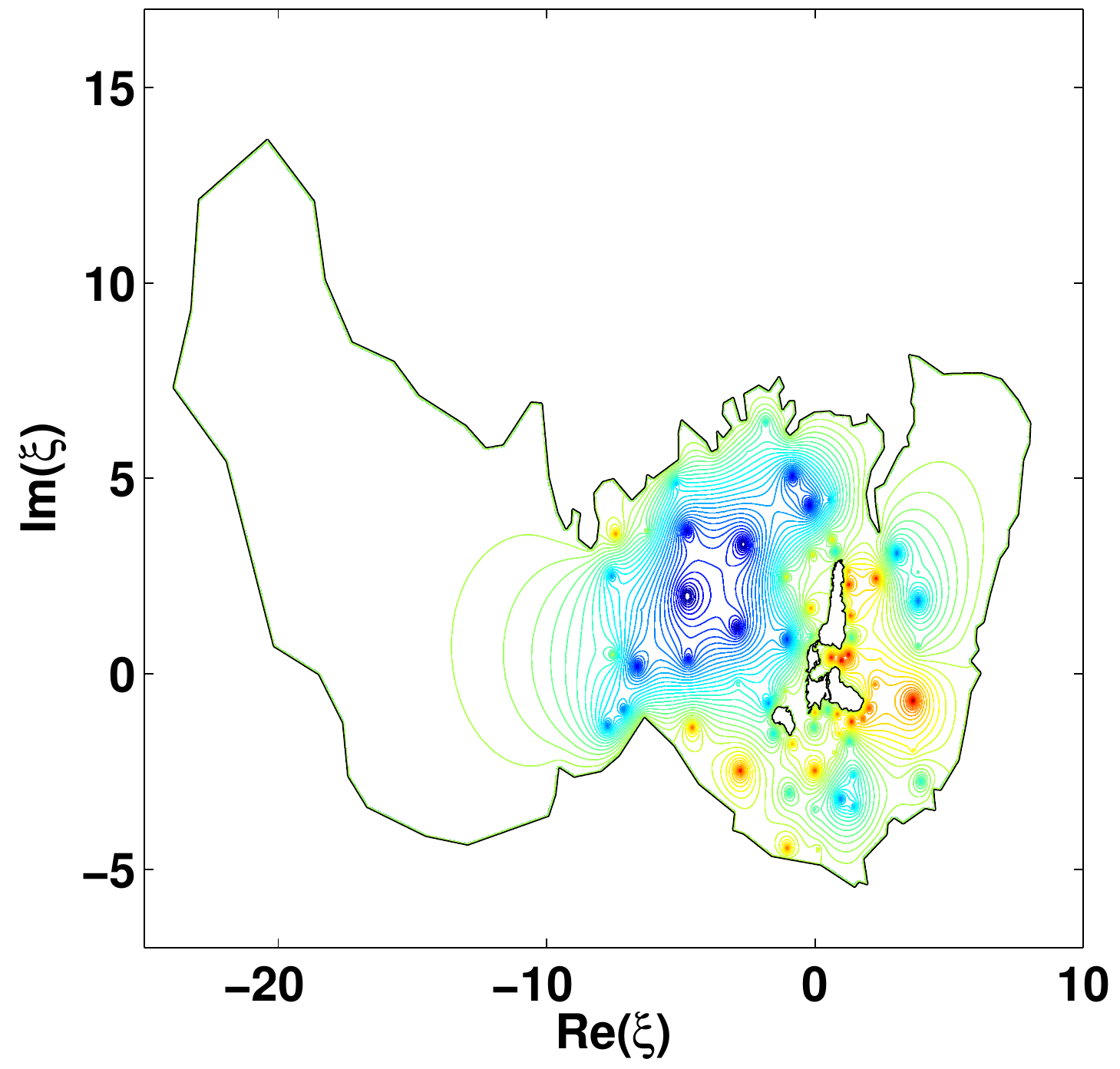} \hspace{.25in}  & \hspace{.25in}
\includegraphics[width=2.in]{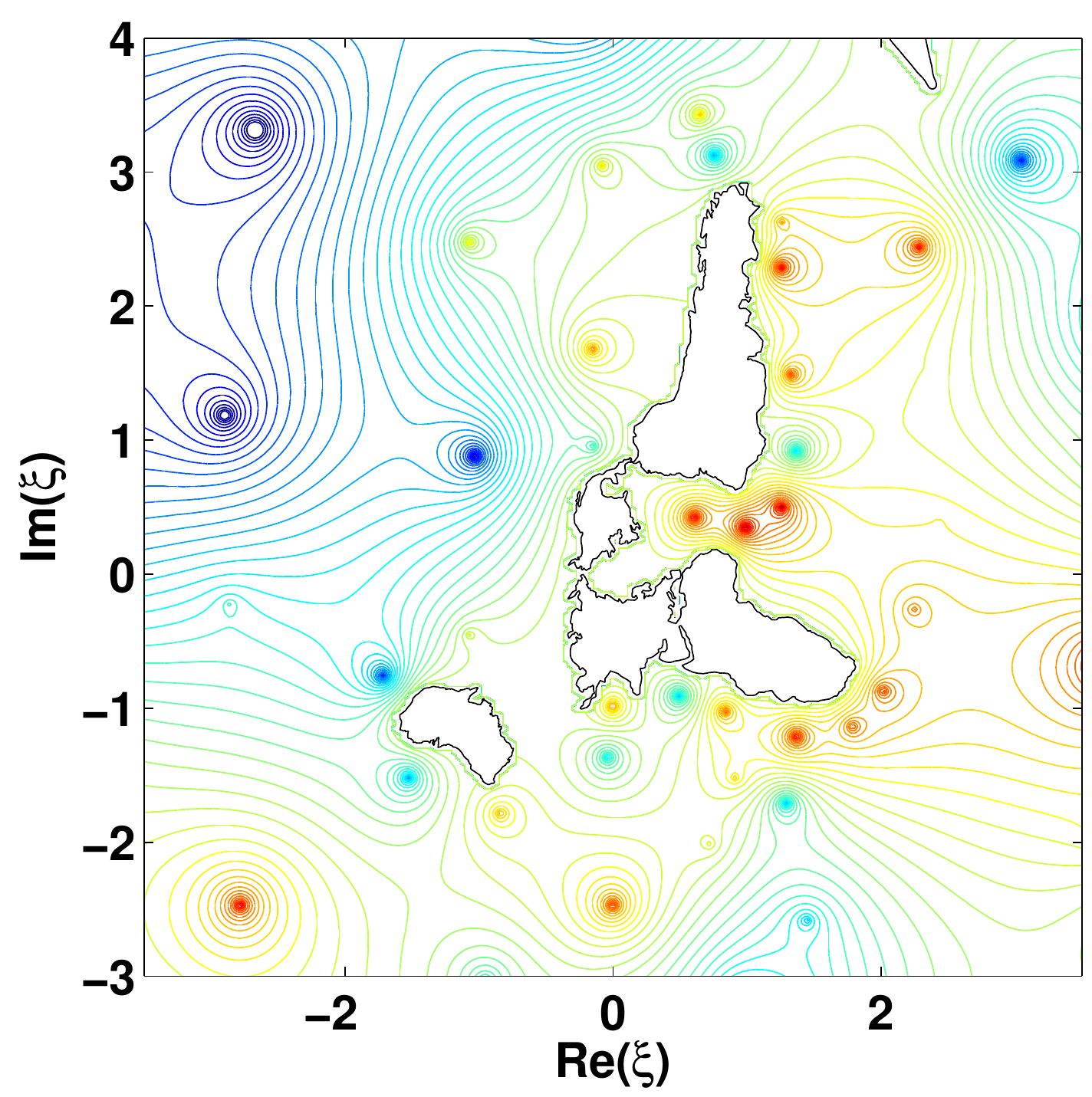}
\end{array}$
\caption{\em The solution to example 2.}
\label{fig5}
\end{figure}

\noindent {\sc Example 2: Point Vortices on Earth.}
To demonstrate the capabilities of our methods on a more realistic geometry, we consider solving\eqr{continuous_sys} in the presence of islands that resemble the Earth's major continents. 
We used coastline data from \textcolor{black}{\tt MATLAB} for Antarctica, the Americas, Australia, Europe and Africa. This data is smoothed and resampled at equally spaced points in arclength (we alter the discretization in \eqr{eqn:disc_sphere} accordingly). The coastline of Antarctica is sampled with 17,736 points, the Americas with 49,728 points, Australia with 23,688 points, and Europe/Africa (which are considered joined) with 50,968 points, for a total number of 142,120 points. For purposes of the stereographic projection, we flip the Earth, so that the ``North Pole'' is in Antarctica (the Earth is flipped back again for plotting purposes). Into the oceans, we place 72 point vortices with strengths randomly generated on the interval $(-2\pi,2\pi)$.
The boundary conditions are set to zero on all of the coastlines.  
To solve this, we first decompose the solution into 
\[
\psi = \psi^* + \sum_{j=1}^{72} \Gamma_j G(\xi,\xi^v_j), 
\]
where $\Gamma_j$ is the strength of point vortex $j$ and $\xi^v_j$ is its location in the stereographic plane.
We then solve for $\psi^*$, which satisfies the Laplace-Beltrami equation, but with non-homogeneous boundary conditions.
GMRES requires 100 iterations and takes 88.8 seconds to converge. The solution is shown in \figr{fig5}.

\noindent \textcolor{black}{\sc Example 3: Performance of FMM. }

Next we demonstrate the importance of using FMM for solving significantly larger problems. \textcolor{black}{An array of $M$ elliptical holes is distributed regularly over the surface of the sphere, with the length of the major and minor axes of these elliptical contours being chosen randomly.} \textcolor{black}{As in Example 1,} the boundary conditions are set according to\eqr{exact_bcs}. \textcolor{black}{We consider situations where the number of islands is} $M=98$, $200$, and $407$; \textcolor{black}{the results} are shown in \tabr{table2}, \tabr{table3}, \tabr{table4}, respectively. \textcolor{black}{ A solution with $M=407$ is shown in \figr{fig3}. }

We should point out, here, a potential drawback of using a 2D FMM on the stereographic plane. As is well known, the FMM introduces a hierarchy of meshes (a quadtree): meshes are adaptively subdivided at increasing levels of refinement until no more than a preset number of discretization points reside in a box at the finest level. \Figr{fig4} shows the quadtree structure for the solution shown in \figr{fig3}, where $N=64$ has been chosen for the purposes of demonstration. Since distances get artificially distorted under the stereographic mapping, the quadtree has to go many levels deep. For $N=64$, 16 levels are needed in the construction of the quadtree; for $N=2048$, that number becomes 22. One might anticipate that introducing a hierarchy of meshes on the unit sphere (an octree) would require far fewer levels to resolve. 
Unfortunately, we do not have a 3D FMM routine that works for the \textcolor{black}{generalized} fundamental solution to the Laplace-Beltrami equation. \textcolor{black}{In future work we will derive and implement a 3D FMM routine suitable for this kernel. }

\begin{table}[h]
\caption{Performance of the algorithm on example 3 with $M=98$. The error is the maximum relative error in the solution sampled at 77 points distributed throughout the sphere,  \# is the number of GMRES iterations required for convergence, CPU Prec indicates the time take for the initial factorization of the preconditioner, CPU Solve is the time needed to solve the linear system. Note that in the largest case with $N=2048$, the total number of unknowns is 200,704.
\label{table2} }
\begin{tabular*}{4in}{@{\extracolsep{\fill}}rrrrl} \hline\noalign{\smallskip}
\multicolumn{1}{c} {$N$} & \multicolumn{1}{c}{\#} & \multicolumn{1}{c}{CPU Prec} 
 & \multicolumn{1}{c}{CPU Solve} &\multicolumn{1}{c}{Error } \\
\noalign{\smallskip}\hline\noalign{\smallskip}
128 & 178 & 0.04 & 36.8 & $4.548\times 10^{-2}$ \\  
256 & 92 & 0.08 & 32.8 & $5.716\times 10^{-4}$ \\  
512 & 92 & 0.2 & 59.1 & $3.995\times 10^{-6}$ \\  
1024 & 92 & 0.33 & 107.4 & $8.744\times 10^{-10}$ \\  
2048 & 92 & 0.65 & 201.5 & $8.924\times 10^{-11}$ \\  
\noalign{\smallskip}\hline
\end{tabular*}
\end{table}
\begin{table}[h]
\caption{Performance of the algorithm on example 3 with $M=200$. Note that in the largest case with $N=2048$, the total number of unknowns is 409,600.
\label{table3} }
\begin{tabular*}{4in}{@{\extracolsep{\fill}}rrrrl} \hline\noalign{\smallskip}
\multicolumn{1}{c} {$N$} & \multicolumn{1}{c}{\#} & \multicolumn{1}{c}{CPU Prec} 
 & \multicolumn{1}{c}{CPU Solve} &\multicolumn{1}{c}{Error } \\
\noalign{\smallskip}\hline\noalign{\smallskip}
128 & 155 & 0.2 & 82.1 & $3.588\times 10^{-2}$ \\  
256 & 86 & 0.3 & 85.3 & $9.973\times 10^{-3}$ \\  
512 & 83 & 0.7 & 160.5 & $4.360\times 10^{-5}$ \\  
1024 & 83 & 1.4 & 312.2 & $6.370\times 10^{-8}$ \\  
2048 & 83 & 2.8 & 622.9 & $7.235\times 10^{-11}$ \\  
\noalign{\smallskip}\hline
\end{tabular*}
\end{table}
\begin{table}[h]
\caption{Performance of the algorithm on example 3 with $M=407$. Note that in the largest case with $N=2048$, the total number of unknowns is 833,536.
\label{table4} }
\begin{tabular*}{4in}{@{\extracolsep{\fill}}rrrrl} \hline\noalign{\smallskip}
\multicolumn{1}{c} {$N$} & \multicolumn{1}{c}{\#} & \multicolumn{1}{c}{CPU Prec} 
 & \multicolumn{1}{c}{CPU Solve} &\multicolumn{1}{c}{Error } \\
\noalign{\smallskip}\hline\noalign{\smallskip}
128 & 125 & 0.8 & 224.5 & $1.727\times 10^{-2}$ \\  
256 & 85 & 1.4 & 301.9 & $1.354\times 10^{-3}$ \\  
512 & 84 & 2.9 & 576.6 & $3.105\times 10^{-5}$ \\  
1024 & 84 & 5.7 & 1132.1 & $2.093\times 10^{-8}$ \\  
2048 & 84 & 11.5 & 2257.0 & $1.459\times 10^{-10}$ \\  
\noalign{\smallskip}\hline
\end{tabular*}
\end{table}

\section{Conclusions}
We have presented \textcolor{black}{an accelerated and preconditioned  integral equation strategy} to solve the Dirichlet boundary value problem for the Laplace-Beltrami equation on a multiply-connected, sub-manifold of the unit sphere. 
\textcolor{black}{By mapping to the stereographic plane, we were able to exploit available fast algorithms for potential problems in $\mathbb{R}^2$ in order to accelerate the solution procedure. }
Our numerical examples confirm that the cost is $O(N)$, where $N$ is the total number of nodes in the discretization \textcolor{black}{\it on each component $\mathcal{C}_k$ of the boundary. We} are able to efficiently compute highly accurate solutions to problems with complicated geometry and of moderate size using only modest computational resources. 
We believe these methods are highly suitable for considering much larger-scale problems of this nature.

We believe there are a number of very fruitful avenues for future investigation, including both theoretical and applied considerations. One immediate future application is to compute the motion of point vortices on the surface of the sphere in the presence of complicated geography. This would extend prior work by other authors \cite{Kidambi,Crowdy2006}, but with more general geometries. 
We also wish to examine solutions on more general manifolds, in which a stereographic projection is either not available, or is not convenient. This would likely require constructing three dimension FMM based on the potentials for the corresponding elliptic operators on these manifolds.   

\noindent
{\bf Appendix A: Continuity of the Kernel.}
Given that the two integral equations\eqr{inteqn:simple} and\eqr{stereo:simple} are equivalent, and that the kernel in\eqr{stereo:simple} is continuous, Lemma 1 follows directly. It is useful, however to have the  expression analogous to\eqr{cont_kern} on the sphere. Here, we derive this expression using elementary arguments.
{\bf Proof: }\\
Let $\s'$ and $\n'$ be the unit tangent and normal vectors at the point $\x'$ lying in the tangent plane of ${\cal S}$. 
We note for future reference that $\x=\er$, $\x \cdot \x = 1$ and that $\x = \n \times \s$. 
Calculating the gradient of $G$ yields
$$
\nabla ' G(\x,\x') = \frac{1}{2\pi} \frac{\x - \x'}{|\x - \x'|^2}, 
$$
which we can decompose into the surface gradient plus a derivative in a direction normal to the sphere. 
Therefore, we can write the surface gradient of the Green's function as 
\[
\nabla' _S G(\x,\x') = \frac{1}{2\pi} \frac{\x - \x'}{|\x - \x'|^2} - (\nabla'  G \cdot \x')  \, \x' , 
\]
and 
\begin{align}
 \frac{\partial \, }{\partial n'} G(\x,\x') & = \frac{1}{2\pi} \frac{\x - \x'}{|\x - \x'|^2} \cdot \n' = \frac{1}{2\pi} \frac{\x - \x'}{|\x - \x'|^2} \cdot \left( \s' \times \x' \right) . \label{kern1}
\end{align}
The above will be shown to be continuous as $\x' \rightarrow \x$ by l'H\^{o}pital's rule. First we note the following identities:
\[
   \s = \frac{d\x}{ds}, \qquad \frac{d\s}{ds} = \kappa \, \N_p .
\]
The second identity is one of the Frenet formulae, where $\kappa$ is the curvature of the curve $C_k$ at the point $\x$ and $\N_p$ is the principal normal to the curve. 
From these, it is straightforward to show that 
\begin{equation}
\frac{d \, }{ds} \left( \s \times \x \right)= \kappa \N_p \times \x . \label{kernel_diag}
\end{equation}
Now, proceeding with the first application of l'H\^{o}pital's rule to evaluate the kernel of the double layer potential at the point of singularity, $\x' = \x$, 
\begin{align*}
   \lim_{\x' \rightarrow \x} \frac{\partial \,}{\partial n'} G(\x,\x')  
      & =  \frac{1}{2\pi} \lim_{\x' \rightarrow \x}  \frac{ -\s' \cdot \left( \s' \times \x' \right) 
                                                             + (\x-\x') \cdot \left( \kappa' \N'_p \times \x'  \right)} 
                                                           { -2 (\x - \x') \cdot \s' } \\ 
                                                        &=  -\frac{1}{4\pi} \lim_{\x' \rightarrow \x}  \frac{  (\x-\x') \cdot \left( \kappa' \N'_p \times \x'  \right)}
                                                           { (\x - \x') \cdot \s' } .
\end{align*}
Proceeding with a second application of l'H\^{o}pital's rule yields
\begin{align*}
   \lim_{\x' \rightarrow \x} \frac{\partial \,}{\partial n'} G(\x,\x')  
      & =  -\frac{1}{4\pi} \lim_{\x' \rightarrow \x}  \frac{  -\s' \cdot 
                                                                 \left( \kappa' \N'_p \times \x'  \right) + O(\x-\x')}
      { -\s' \cdot \s' + O(\x-\x')}  \\
  &= -\frac{1}{4\pi} \s \cdot  \left( \kappa \N_p \times \x  \right) 
\end{align*}


\bibliographystyle{spmpsci}

\end{document}